\documentclass[12pt]{article}
\usepackage{pictexwd,dcpic}
\usepackage[all,cmtip]{xy}
\usepackage{stmaryrd}
\usepackage{amsfonts}
\usepackage{mathrsfs}
\usepackage{amsmath}
\usepackage{amssymb,amsmath}
\begin{document}

\setcounter{equation}{0}

\newcommand{\wt}{\mbox{wt}}
\newcommand{\spa}{\mbox{span}}
\newcommand{\Res}{\mbox{Res}}
\newcommand{\End}{\mbox{End}}
\newcommand{\Ind}{\mbox{Ind}}
\newcommand{\Hom}{\mbox{Hom}}
\newcommand{\Mod}{\mbox{Mod}}
\renewcommand{\theequation}{\thesection.\arabic{equation}}

\def \Aut{{\rm Aut}}
\def \Z{\Bbb Z}
\def \M{\Bbb M}
\def \C{\Bbb C}
\def \R{\Bbb R}
\def \Q{\Bbb Q}
\def \N{\Bbb N}
\def \ann{{\rm Ann}}
\def \<{\langle}
\def \o{\omega}
\def \O{\Omega}
\def \M{{\cal M}}
\def \1t{\frac{1}{T}}
\def \>{\rangle}
\def \t{\tau }
\def \a{\alpha }
\def \e{\epsilon }
\def \l{\lambda }
\def \L{\Lambda }
\def \g{\gamma}
\def \b{\beta }
\def \om{\omega }
\def \o{\omega }
\def \cg{\chi_g}
\def \ag{\alpha_g}
\def \ah{\alpha_h}
\def \ph{\psi_h}
\def \nor{\vartriangleleft}
\def \V{V^{\natural}}
\def \voa{vertex operator algebra\ }
\def \v{vertex operator algebra\ }
\def \1{{\bf 1}}
\def \be{\begin{equation}\label}
\def \ee{\end{equation}}
\def \pf {\noindent {\bf Proof:} \,}
\def \bl{\begin{lem}\label}
\def \el{\end{lem}}
\def \ba{\begin{array}}
\def \ea{\end{array}}
\def \bt{\begin{thm}\label}
\def \et{\end{thm}}
\def \ch{{\rm ch}}
\def \br{\begin{rem}\label}
\def \er{\end{rem}}
\def \ed{\end{de}}
\def \bp{\begin{prop}\label}
\def \ep{\end{prop}}
\def\QED{\hfill$\Box$}

\newtheorem{th*}{Theorem}
\newtheorem{conj}[th*]{Conjecture}
\newtheorem{ree}[th*]{Remark}
\newtheorem{thm}{Theorem}[section]
\newtheorem{prop}[thm]{Proposition}
\newtheorem{coro}[thm]{Corollary}
\newtheorem{lem}[thm]{Lemma}
\newtheorem{rem}[thm]{Remark}
\newtheorem{de}[thm]{Definition}
\newtheorem{hy}[thm]{Hypothesis}

\begin{center}

{\Large {\bf On rationality of vertex operator superalgebras}} \\
\vspace{0.5cm}

Chongying Dong\footnote{Supported by NSF grants and a faculty
research fund from the University of California at Santa Cruz.}
\\
Department of Mathematics\\ University of
California, Santa Cruz, CA 95064 \\
Jianzhi Han\\
Department of Mathematics, Sichuan University\\ Chengdu 610064,
China
\end{center}
\hspace{1.5 cm}
\begin{abstract}
In this paper we prove that
 $g$-rationality of  a vertex operator
superalgebra $V=V_{\bar0}\oplus V_{\bar1}$ for all $g\in G$ imply
rationality of $V^G$, and also imply that each irreducible
$V^G$-module is a submodule of an irreducible $g$-twisted
$V$-module for some $g\in G$,  where $G$ is any finite abelian
subgroup of $Aut(V)$. We also prove that for any finite solvable
$G$, rationality of $V^G$ implies $g$-rationality of $V$ for any
$g\in G$.
\end{abstract}

 \section{Introduction}
 \setcounter{section}{1}\setcounter{equation}{0}

This paper is a  continuation of \cite{DH}, studying the rationality of a vertex operator superalgebra $V=V_{\bar0}\oplus V_{\bar1}$. In \cite{DH}, we conjectured that the rationality of $V$, the
$\sigma$-rationality of $V$ and the rationality of $V_{\bar0}$
are equivalent to each other, where $\sigma$ is the canonical
automorphism of $V$ defined by $\sigma|_{V_{\bar i}}=(-1)^i$. We prove in this paper that the rationality of
$V$ together with the $\sigma$-rationality of $V$ is equivalent to the rationality of $V_{\bar0}.$
In fact, we prove a stronger result concerning the general orbifold
theory for a vertex operator superalgebra $V:$
\begin{th*}\label{t:I}
Let $V$ be a simple vertex operator superalgebra of CFT type and
$G$ a finite abelian subgroup of Aut$(V)$. Suppose that $V^G$ is
$C_2$-cofinite. Then $V^G$ is rational if and only if $V$ is
$g$-rational for all $g\in G$.
\end{th*}

Proving the main result
depends heavily on the decomposition of
$$V=\oplus_{\chi\in{\rm
Irr}G}V^\chi$$ for the  abelian group $G$ (see \cite{DLM1}). Here,
$V^{\chi}$ is the sum of irreducible $G$-submodules of $V$ which
afford the character $\chi.$ A $V$-module $M=\oplus_{\chi\in{\rm
Irr}G} M^\chi$ satisfying  certain conditions was considered in
\cite{DLM3} and proved to be completely reducible under the
assumption that each $V^\chi$ is a simple current.  Using the
decomposition $V=\oplus_{\chi\in{\rm Irr}G}V^\chi$  and the same
idea as in \cite{DLM3} we  prove that for any finite abelian
group $G$,   $V$ is $g$-rational for any $g\in G$ if $V^G$ is
rational. Applying   this result again we  obtain
that the above conclusion  holds true for any finite solvable group
$G$. Note that in our case we do not assume that $V^G$ is
$C_2$-cofinite. So even $V^G$ is rational we cannot claim that
each $V^{\chi}$ is a simple current \cite{DJX} which will make the
proof much easier.

The proof of the opposite direction involves another well known
conjecture in the orbifold theory \cite{DVVV}: Each irreducible
$V^G$-module is a submodule of an irreducible $g$-twisted
$V$-module for some $g\in G$. We establish the conjecture for
finite abelain groups in this paper:
\begin{th*}\label{conjec}
Let $V$ be a simple vertex operator superalgebra of CFT type and
$G$ a finite abelian subgroup of Aut$(V)$. Suppose that $V^G$ is
$C_2$-cofinite and $V$ is $g$-rational for any $g\in G$.  Then
each irreducible $V^G$-module is a $V^G$-submodule of  an irreducible
$g$-twisted $V$-module for some $g\in G$.
\end{th*}

This implies that one can construct a (twisted) $V$-module from
any given module for $V^G$, which is the key in our argument of
rationality of $V^G$ from $g$-rationality of $V$ for any $g\in G.$
The  proof uses the generalized tensor product theory developed in
\cite{HLZ}.

As an application, we give a vertex operator superalgebra
theoretical proof of the $SL(\!2,\Z\!)$-invariance of the partition
functions  of type $D_{2\rho +1}$ given in \cite{CIZ}. The other
types except $E_7$ have vertex operator algebra or conformal nets
interpretation already \cite{KaL}. Motivated by this example we
make the following conjecture on  vertex operator superalgebras
\begin{conj}\label{conje-AAA}
Let $V$ be a simple vertex operator superalgebra which is
rational, $\sigma$-rational, $C_2$-cofinite and self-dual, i.e.,
$V^\prime\cong V$. Suppose that $M_i$ for $i=1,...,n$ and
$W_j$ for $j=1,...,m$ are all the irreducible
$V$-modules and irreducible $\sigma$-twisted $V$-modules,
respectively.
Then$$\sum_{i=1}^n\big(|Tr_{M_{i}}q^{L(0)-\frac{c}{24}}|^2+|Tr_{M_{i}}\sigma
q^{L(0)-\frac{c}{24}}|^2\big)+\sum_{j=1}^m|Tr_{W_{i}}q^{L(0)-\frac{c}{24}}|^2$$
is $SL(2,\Z)$-invariant.
\end{conj}

Note from \cite{DZ2} that the space spanned by $Tr_{M_{i}}q^{L(0)-\frac{c}{24}},$
$Tr_{M_{i}}\sigma q^{L(0)-\frac{c}{24}},$ $ Tr_{W_{j}}q^{L(0)-\frac{c}{24}}$ for all $i,j$ affords a representation of $SL(2,\Z).$ As in \cite{DLN}, the conjecture simply says that this representation
of $SL(2,\Z)$ is unitary.

This paper will be organized as follows: In Section 2, we shall
review some basic definitions and notation. Theorem 2 is proved in Section 3 and Theorem 1 is established in Section
4. We apply the general result on vertex operator superalgebras obtained in this paper to give an interpretation of the affine partition functions on type $D_{2\rho+1}$ proposed in \cite{CIZ} in Section 5.

\section{
Basics}\setcounter{section}{2}\setcounter{thm}{0}\setcounter{equation}{0}

\bigskip
In this section we review the $P(z)$-tensor products of  certain
generalized (twisted) $V$-modules (cf. \cite{HLZ1}) and related
results.

We assume that the reader is familiar with various notions of
twisted modules (see \cite{DLM}, \cite{DZ1}, \cite{DH}). Let $V$
be a vertex operator superalgebra and $g$ an automorphism of $V$.
\begin{de}
A weak $g$-twisted $V$-module $M$ is said to be of length $n$ if
there exist submodules $0=M^0\subset M^1\subset\cdots\subset
M^n=M$ such that $M^{i+1}/M^i$ for $i=1,2,\cdots, n-1$ are
irreducible $g$-twisted $V$-modules.
\end{de}

Note that a weak $g$-twisted $V$-module $M$ of finite length, in fact, is admissible. There exist $\lambda_1,...,\lambda_s\in\C$ such that
$$M=\bigoplus_{i=1}^s\bigoplus_{n\geq 0}M_{\lambda_i+n}$$
where $M_{\lambda_i+n}=\{w\in M|(L(0)-\lambda_i-n)^kw=0, {\rm for\
some}\ k\geq 0\}$ is the generalized eigenspace of $L(0)$ with
eigenvalue $\lambda_i+n.$  Moreover, $M_{\lambda_i+n}$ is
necessarily finite dimensional. That is, $M$ is $\C$-graded with
finite dimensional homogeneous subspaces.

Let $\mathcal F_{g,V}$ denote the category of weak $g$-twisted
$V$-modules of finite length. We simply write $\mathcal F_{g, V}$
as $\mathcal F_V$ if $g=1$.

Let $g_i$ be mutually commutative automorphisms of $V$ with finite
order $T_i$ for $i=1,2,3$. Set $V^{(s,t)}=\{v\in V\,|\,
g_1v=e^{\frac{2\pi s i}{T_1}}v\,\, {\rm and}\,\,
g_2v=e^{\frac{2\pi t i}{T_2}}v\}$
\begin{de}
Let  $W_i$ be in  $\mathcal F_{g_i,V}$ for $i=1,2,3$ and  $z$  a nonzero complex number.  A
$P(z)$-intertwining map of type $\left(\begin{array}{cc}
    W_3 \\
  W_1\ W_2 \\
\end{array}\right)$  is a linear map $F:
W_1\otimes W_2\rightarrow \overline W_3$ such that
\begin{eqnarray*}&&x_0^{-1}\delta\big(\frac{x_1-z}{x_0}\big)^{\frac{s}{T_1}}x_0^{-1}\delta\big(\frac{x_1-z}{x_0}\big)Y_3(v,x_1)F(w_1\otimes
w_2)\\
&&\ \ \ \ \ \ \ \ \ -(-1)^{ij}x_0^{-1}\delta\big(\frac{z-x_1}{-x_0}\big)^{\frac{s}{T_1}}x_0^{-1}\delta\big(\frac{z-x_1}{-x_0}\big)F(w_1\otimes
Y_2(v,x_1)w_2)\\
&&\ \ \ \ \ \ \ \ \
=z^{-1}\delta\big(\frac{x_1-x_0}{z}\big)^{-\frac{t}{T_2}}z^{-1}\delta\big(\frac{x_1-x_0}{z}\big)F(Y_1(v,x_0)w_1\otimes
w_2)\end{eqnarray*} for any $v\in V_{\bar i}\cap V^{(s,t)}$,
$w_1\in (W_1)_{\bar j}$ and $w_2\in W_2$, where $\overline W_3$ is
the completion $\prod_{n\in\C}(W_3)_n$ of $W_3$ with respect to
the $\C$-grading.
\end{de}

For any $W_i$  in $\mathcal F_{g_i,V}$, a
$P(z)$-product of $W_1$ and $W_2$ is a $g_1g_2$-twisted $V$-module $(W_3, Y_3)$ for
which $W_3\in\mathcal F_{g_1g_2,V}$ equipped with a $P(z)$-intertwining map
$F$ of type $\left(\begin{array}{cc}
    W_3 \\
  W_1\ W_2 \\
\end{array}\right)$  which  is denoted  by
$(W_3, Y_3; F)$. Let $(W_4, Y_4; H)$ be another $P(z)$-product of
$W_1$ and $W_2$. A {\em morphism} from $(W_3, Y_3; F)$ to $(W_4,
Y_4, H)$ is a module map $\eta$ from $W_3$ to $W_4$ such that
\begin{equation}
H=\overline{\eta}\circ F,
\end{equation}
where $\overline{\eta}$ is the natural map from $\overline W_3$ to
$\overline W_4$ uniquely extending $\eta$.

\begin{de}
A $P(z)$-tensor product of $W_1$ and $W_2$ is a $P(z)$-product
$$(W_1\boxtimes_{P(z)}W_2, Y_{P(z)}; \boxtimes_{P(z)})$$ such that
for any $P(z)$-product $(W_3, Y_3; F)$, there is a unique morphism
from $(W_1\boxtimes_{P(z)}W_2, Y_{P(z)}; \boxtimes_{P(z)})$ to
$(W_3, Y_3; F)$.  The $V$-module $(W_1\boxtimes_{P(z)}W_2,
Y_{P(z)})$ is called a $P(z)$-tensor product module of $W_1$ and
$W_2$.
\end{de}

\begin{rem}
Though the the $P(z)$-tensor product $W_1\boxtimes_{P(z)}W_2$ and
the usual box tensor product $W_1\boxtimes W_2$ seem different, as
a matter of fact, they are the same thing provided that both of
them exist (cf. \cite{Li0}).
\end{rem}

It is clear from the definition that the $P(z)$-intertwining maps are
similar to the logarithmic intertwining operators (cf. \cite{HLZ1}).
In fact, we can establish a one-to-one correspondence between  logarithmic intertwining operators  and $P(z)$-intertwining maps of the same type. To be more explicit, let $\mathcal Y$ be a logarithmic
intertwining operator of type $\left(\begin{array}{cc}
    W_3 \\
  W_1\ W_2 \\
\end{array}\right)$. Define a linear map $F_{\mathcal Y,p}^{P(z)}: W_1\otimes W_2\rightarrow \overline W_3$
by
\begin{eqnarray*}
F_{\mathcal Y, p}^{P(z)}(w_1\otimes w_2)=\mathcal
Y(w_1,e^{l_p(z)})w_2
\end{eqnarray*} where $l_p(z)={\rm log}z+2p\pi i$. Here and below we always choose log$z$ to be
log$|z|+i{\rm arg}z$ with $0\leq {\rm arg}z< 2\pi$. Conversely,
let $F$ be a given $P(z)$-intertwining map, $w_1\in W_1$, $w_2\in
W_2$ be homogenous and $n\in\C$. We define a linear map $\mathcal
Y_{F, p}^{P(z)}$ from $W_1\otimes W_2$ to $W_3[{\rm log}x]\{x\}$,
which is given by
\begin{eqnarray*}
\mathcal Y_{F,
p}^{P(z)}(w_1,x)w_2=y^{L(0)}x^{L(0)}F(y^{-L(0)}x^{-L(0)}w_1\otimes
y^{-L(0)}x^{-L(0)}w_2)|_{y=e^{-l_p(z)}}
\end{eqnarray*}for any $w_1\in W_1$ and $w_2\in W_2$.

The following result is proved in \cite{HLZ}:
\begin{prop}\label{interwing-isomorphic}
For $p\in\Z$, the correspondence $\mathcal Y\mapsto F_{\mathcal Y,
p}^{P(z)}$ is a linear isomorphism from the vector space $\mathcal
V_{W_1,W_2}^{W_3}$ of logarithmic intertwining operators of type
$\left(\begin{array}{cc}
    W_3 \\
  W_1\ W_2 \\
\end{array}\right)$  to the vector space
$\mathcal M[P(z)]_{W_1,W_2}^{W_3}$ of $P(z)$-intertwining maps of
the same type. Its inverse map is given by $F\mapsto \mathcal
Y_{F, p}^{P(z)}$.
\end{prop}

From now on  we identify $P(z)$-intertwining maps with their corresponding
logarithmic intertwining operators  if necessary.

For a positive integer $n$ and a weak $V$-module $M$, we say that $M$
is {\em $C_n$-cofinite} (cf. \cite{Z}
and \cite{Li4}) if dim$V/C_n(M)<\infty$, where
\begin{equation*}
C_n(M)=\langle u_{-n}w|u\in\coprod_{i>0}V_i, w\in M\rangle.
\end{equation*}
In particular, we can define $C_n$-cofinite vertex operator algebra $V$ for $n>1.$ But for $n=1$ we have to modify $C_1(V)$ as follows
 $$C_1(V)=\langle u_{-1}v, L(-1)w|u,v\in \coprod_{i>0}V_i, w\in V\rangle,$$
otherwise, $C_1(V)$ could be the whole $V.$ This also gives a $C_1$-cofinite vertex operator algebra.

\begin{rem} Let $G$ be a finite subgroup of $Aut(V)$.  Assume that $V^G$ is $C_2$-cofinite and $V$ is of  CFT type, i.e.,
$V=\coprod_{n\geq0}V_n$ and $V_0=\C\bf1$.   Then the $P(z)$-tensor
product $W_1\boxtimes_{P(z)} W_2$ exists for $W_i\in \mathcal
F_{g_i,V}$ and $z\in\C^*$, which  lies in $\mathcal F_{g_1g_2,V}$
(see \cite{H}). Moreover,  the $P(z)$-tensor products satisfy the
associative property \cite{HLZ}: For any $z_1$, $z_2$ complex
numbers satisfying
$$|z_1| > |z_2| > |z_1- z_2| > 0,$$  there exists a unique
natural isomorphism
$$\mathcal A_{P(z_1), P(z_2)}^{P(z_1-z_2),P(z_2)}: W_1\boxtimes_{P(z_1)}(W_2\boxtimes_{P(z_2)} W_3)\rightarrow
(W_1\boxtimes_{P(z_1-z_2)}W_2)\boxtimes_{P(z_2)}W_3$$ such that
$\mathcal A_{P(z_1),
P(z_2)}^{P(z_1-z_2),P(z_2)}\big(w_1\boxtimes_{P(z_1)}(w_2\boxtimes_{P(z_2)}
w_3)\big)=(w_1\boxtimes_{P(z_1-z_2)}w_2)\boxtimes_{P(z_2)}w_3$ for
$w_i\in W_i$.
\end{rem}

\section{Irreducible $V^G$-modules}
\setcounter{section}{3}\setcounter{thm}{0}\setcounter{equation}{0}
\bigskip

We prove Theorem 2 in this section. That is, any irreducible
$V^G$-module is a submodule of some $g$-twisted $V$-module.

Let $V$ be a simple vertex operator superalgebra and $G$ a finite
abelian subgroup of $Aut(V)$. Then by \cite{DLM1} we have the
decomposition $V=\oplus_{\chi\in {\rm Irr}G}V^\chi$, where each
$V^\chi=\{v\in V|gv=\chi(g)v\}$ is an irreducible $V^G$-module on
which $G$ acts according to the character $\chi$.

The following elementary lemma  is the key in  proof of several
results.
\begin{lem}\label{lemma-1}
Let $M$ be an admissible $g$-twisted  $V$-module and $L$  an
irreducible $V^G$-submodule of $M$. Then $V^\chi\cdot L=\langle
u_nw\ |\ u\in V^\chi, w\in L, n\in \Q\rangle$ is irreducible as
$V^G$-module for any $\chi\in$ {\rm Irr}$G$.
\end{lem}
\noindent{\em Proof.}\ \ Clearly, $V^\chi\cdot L$ is an admissible
$V^G$-module. To show the irreducibility of $V^\chi\cdot L$, let
$T$ be a nonzero $V^G$-submodule of $V^\chi\cdot L$, we must prove
$T=V^\chi\cdot L$. Notice that
$$0\neq V^{\chi^{-1}}\cdot T\subset V^{\chi^{-1}}\cdot(V^\chi\cdot L)=(V^{\chi^{-1}}\cdot V^\chi)\cdot
L=V^G\cdot L=L.$$ Then $V^{\chi^{-1}}\cdot T=L$, as $L$ is
irreducible. Applying the action of $V^\chi$ to this equality we
have $T=V^G\cdot T=V^\chi\cdot(V^{\chi^{-1}}\cdot T)=V^\chi\cdot
L$, as desired. \QED
\bigskip

Recall that an admissible $g$-twisted $V$-module $W$  is called
{\em projective} if for any admissible $g$-twisted $V$-module $M$
and any epimorphism $\pi: M\rightarrow W$, there exists a
homomorphism $\phi: W\rightarrow M$ such that $\pi\phi={\rm
Id}_W$. In the following lemma we are going to use the generalized
Zhu's algebra $A_{g, n}(V)$ \cite{DLM4} to show the existence of
projective admissible $g$-twisted $V$-modules.

\begin{lem}\label{lemma-projective} Let $V$ be a vertex operator superalgebra of CFT type and $G$  any finite
subgroup of $Aut(V)$.  Assume that $V^G$ is $C_2$-cofinite. Then
each object in $\mathcal F_{g, V}$ is  a homomorphic image of a
projective object in the same category.
\end{lem}
\noindent{\em Proof.} \ \ For any projective $A_{g, n}(V)$-module
$U$, we assert that the admissible $g$-twisted $V$-module
$\overline M_{n}(U)$ constructed in Section 4 of \cite{DLM4} is
projective. For any admissible $g$-twisted $V$-module $M$ and any epimorphism $f: M\rightarrow
\overline M_n(U)$, we get an $A_{g, n}(V)$-epimorphism from $\Omega_n(M)$ to $\Omega_n(\overline M_{n}(U))$
where
$$\Omega_n(W)=\{w\in W\,|\, u_kw=0\ {\rm for \
homogenous }\ u\in V, {\rm wt}u-k-1<-n\}$$ for any weak
$g$-twisted $V$-module $W.$ Note that $U$ is an $A_{g,
n}(V)$-submodule of $\Omega_n(\overline M_{n}(U)).$ As a result
there exists an $A_{g, n}(V)$-submodule $\tilde U$ of
$\Omega_n(M)$ such that $f|_{\tilde U}: \tilde U\to U$ is an
$A_{g, n}(V)$-epimorphism. Since $U$ is projective, there exists
an $A_{g, n}(V)$-homomorphism $\tilde h: U\rightarrow \tilde U$
such that $f\circ \tilde h={\rm Id}_U$. Now it follows from the
universal property of $\overline M_n(U)$ (see Theorem 4.1 of
\cite{DLM4}) that the $A_{g,n}(V)$-homomorphism  $\tilde h$ can be extended
to a $V$-homomorphism $h: \overline M_n(U)\rightarrow M$ such that
$f\circ h={\rm Id}_{\overline M_n(U)}$. That is, $\overline
M_n(U)$ is projective, proving our assertion.

As $V^G$ is $C_2$-cofinite,  then so is  $V^{\langle g\rangle }$ \cite{ABD}.
In particular, $V^{\langle g\rangle }$ satisfies the condition in
Corollary 3.17 of \cite{H}.  If $U$ is finite dimensional, then by
Corollary 3.17 of \cite{H} we obtain $\overline M_n(U)\in \mathcal
F_{V^{\langle g\rangle}}$, and hence each generalized eigenspace
$(\overline M_n(U))_\mu$ of $L(0)$ with eigenvalue weight $\mu$ is
finite dimensional. Then it follows from the proof of Proposition
9.4 of \cite{DLM} that $\overline M_n(U)\in \mathcal F_{g,V}$.

Let $ W=\oplus_{n\in \frac{1}{T}\Z_+} W(n) $ be in $\mathcal F_{g,
V}$, where $g\in G$ and $T$ is the order of $g$. It follows from
the proof of Theorem 3.3 of \cite{DH} that there exists a positive
integer $N$ depending only on $g$ such that $W$ is generated by
$\oplus_{n=0}^{{\frac N T}} W(n)$. For each $n$ we can choose a
 finite dimensional projective $A_{g,n}(V)$-module $U_n$ with a
homomorphic image $W(n )$. It follows from the previous arguments
that $\oplus_{n=0}^{\frac{N}{T}}\overline M_n(U_n)$ is a
projective
 object in $\mathcal F_{g, V}$ with a homomorphic image
$ W$. This completes the proof. \QED

\bigskip

There is no induction theory in general for vertex operator
superalgebras. The situation is even worse when  twisted
modules are involved. The following result using the $P(z)$-tensor
product gives some kind of induction which will be used later.

\begin{thm}\label{thm-submodule}
Let $V$ be a simple vertex operator superalgebra of CFT type and
$G$ a finite abelian subgroup of Aut$(V)$. Suppose that $V^G$ is
$C_2$-cofinite. Then every indecomposable projective   object in
${\mathcal F}_{V^G}$ is a $V^G$-homomorphic image of  a sum of
some admissible $g$-twisted $V$-module  for  $g\in G$.
\end{thm}

\noindent{\em Proof.} First consider a special case, that is, $G$
is a cyclic group generated by $g$. The basic ideal of the proof
for this case comes from \cite{M0}. Let $p$ be the order of $g$.
Set $V^{j}=\{v\in V| gv=e^{\frac{2j\pi i}{p}}v\}$ for $0\leq j\leq
p-1$.  Let $W$ be an indecomposable projective  $V^G$-module of
finite length. Note that we have the following two
$V^G$-epimorphisms
$$\underbrace{V^1\boxtimes V^1\boxtimes\cdots \boxtimes
V^1}_p\rightarrow V^G$$ and
\begin{equation}\label{proj-epi1}\underbrace{V^1\boxtimes V^1\boxtimes\cdots
\boxtimes V^1}_p\boxtimes W\rightarrow W.\end{equation} Since $W$
is projective,  $W$ is isomorphic to a direct summand (which is
also a $V^G$-submodule module) of $\underbrace{V^1\boxtimes
V^1\boxtimes\cdots \boxtimes V^1}_p\boxtimes W$ and  we may assume
that $W$ is a $V^G$-submodule of $\underbrace{V^1\boxtimes
V^1\boxtimes\cdots \boxtimes V^1}_p\boxtimes W$. For $0\leq i\leq
p$, let $Q_i$ be an indecomposable  direct summand of
$\underbrace{V^1\boxtimes V^1\boxtimes\cdots \boxtimes
V^1}_i\boxtimes W$ such that $Q_{i+1}$ is a direct summand of
$V^1\boxtimes Q_i$  and $Q_P=W$. So we obtain another epimorphism:
\begin{equation}\label{proj-epi2}
\underbrace{V^1\boxtimes V^1\boxtimes\cdots \boxtimes
V^1}_p\boxtimes W\rightarrow \underbrace{V^1\boxtimes
V^1\boxtimes\cdots \boxtimes V^1}_{p-1}\boxtimes
Q_1\rightarrow\cdots\rightarrow V^1\boxtimes Q_{p-1}\rightarrow W.
\end{equation}

 Set $$Q=\bigoplus_{i=0}^{p-1}Q_i.$$  Let $\mathcal {Y}^i$ be  the natural logarithmic intertwining
operator of type $\left(
\begin{array}{cc}
          Q_{i+1} \\
         V^1\  \ \ Q_i
                          \end{array} \right)$ for $0\leq i\leq p-1$.
Then it follows from the homomorphisms (\ref{proj-epi1}) and
(\ref{proj-epi2}) that
\begin{eqnarray}\label{eqnarray-asso}
&&\langle w^\prime, \mathcal
Y^{p-1}(u^{p-1},z_{p-1})\cdots\mathcal
Y^{0}(u^{0},z_{0})w\rangle\nonumber\\
\label{commua-dia}&=&\langle w^\prime, Y(Y(\cdots Y(u^{p-1},
z_{p-1}-z_{p-2})u^{p-2}, \cdots),z_{1}-z_{0})u^0, z_0) w\rangle.
\end{eqnarray} holds true  for all
$w^\prime \in W^\prime$, $u^i\in V^1$ and $w\in W$ in the region
$|z_{p-1}|>|z_{p-2}|>\cdots>|z_0|>|z_1-z_0|>|z_2-z_1|>\cdots>|z_{p-1}-z_{p-2}|>0$.
Since $V$ is of CFT type and $V^G$ is $C_2$-cofinite, then by
Corollary 3.12 of \cite{KL} and Proposition 2.9 of \cite{H} every
object in $\mathcal F_{V^G}$ satisfies the $C_1$-cofiniteness  and
its weights are all rational numbers by Corollary 5.10 of
\cite{M1}. Hence the left side of (\ref{eqnarray-asso}) satisfies
some  differential equations of regular singular points, and as a
result,  by the theory of differential equations with regular
singular points, can be analytically extended to a multivalued
analytic function
in the region given by $z_1\neq 0, z_2\neq 0, z_1\neq z_2$ with
regular singular points at $z_i=0$ and $z_i=z_j$ (see \cite{H2},
\cite{HLZ} and Appendix B of \cite{K}).

Next we shall define the vertex operator $Y_Q(\cdot,z)$ of $V$ on
the entire $Q$. First, the action of $V^G\oplus V^1$ on $Q$ can be
defined canonically. For $i\geq 2$ and $q_j\in Q_j$, define
$$Y_Q(Y(\cdots (Y(u^1,z_1)u^2,z_2),\cdots,z_{i-1})u^i,z_i)q_j\in
Q_{j+i}$$ (here and below the index $i+j$ is understood as $i+j\
({\rm mod}\ p)$) as follows
\begin{eqnarray}&&\langle q, \mathcal Y^{i+j}(u^1,z_1)\mathcal
Y^{i+j-1}(u^2,z_2)\cdots\mathcal Y^j(u^i,z_i)q_j\rangle\nonumber\\
\label{eq-de-action} &=&\langle q,
Y_Q(Y(\cdots(Y(u^1,z_1-z_2)u^2,z_2-z_3),\cdots,z_{i-1}-z_i)u^i,z_i)q_j\rangle.
\end{eqnarray}
Since all the coefficients of $Y(Y(\cdots
(Y(u^1,z_1)u^2,z_2),\cdots,z_{i-1})u^i,z_i)$ for
$u^j\in V^1$  span $V^i$, as a matter of fact, we have already defined the
action of $V^i$ on $Q$.
 In particular,
we have
\begin{equation}\label{skew-sy}
\langle q^\prime, \mathcal Y^{i+1}(v,z_2)\mathcal
Y^i(u,z_1)q_i\rangle=\langle q^\prime, Y_Q(Y(v,z_2-z_1)u,
z_1)q_i\rangle.
\end{equation}for any $u,v\in V^1$, $q_i\in Q_i$ and $q^\prime \in
Q$. Applying the skew symmetry to the right hand side of (\ref{skew-sy})  and  using the fact that there exists a positive integer $n$ such
that
\begin{equation*}
(z_1-z_2)^nY(u, z_1-z_2)v=(z_1-z_2)^nY(u, -z_2+z_1)v,
\end{equation*}
one has
\begin{equation}\label{weak-comm1}
\langle q^\prime, \mathcal Y^{i+1}(v,z_2)\mathcal
Y^i(u,z_1)q_i\rangle=\langle q^\prime, \mathcal
Y^{i+1}(u,z_1)\mathcal Y^i(v,z_2)q_i\rangle.
\end{equation} So far we have already obtained the
vertex operator $Y_Q(\cdot,z)$. Indeed, the weak associativity
holds for $Y_Q(\cdot,z)$. To see this, by Theorem 3.3.10 of
\cite{X} it is enough to show
\begin{equation}\label{eq-ass}\langle q^\prime, Y_Q(v,z_0+z_i)Y_Q(u,z_i)q\rangle
=\langle q^\prime, Y_Q(Y(v,z_0)u,z_i)q\rangle
\end{equation} for any $v\in V^1$,
$u\in V^i $, $q\in Q$ and $q^\prime\in Q^\prime$. We define the
adjoint vertex operators $Y^\prime_Q(a,z)$ \cite{FHL} on $Q^\prime$ so that
\begin{equation*}
\langle Y^\prime_Q(a,z)q^\prime,  q\rangle=\langle q^\prime,
Y_Q(e^{zL(1)}(-z^{-2})^{L(0)}a,z^{-1})q\rangle
\end{equation*}for $a\in V$, $q\in Q$ and $q^\prime\in Q^\prime$. It follows immediately  (see Section 3 of \cite{HL1}) that
\begin{equation}\label{eq-adjoint}
\langle Y^\prime_Q(e^{zL(1)}(-z^{-2})^{L(0)}a,z^{-1})q^\prime,
q\rangle=\langle q^\prime, Y(a,z)q\rangle.
\end{equation}Here we do not require  $W$ to carry any
structure of a $V$-module.

 Assume
$$u={\rm Res}_{z_1}\cdots {\rm Res}_{z_{i-1}}z_1^{k_1}\cdots
z^{k_{i-1}}_{i-1}Y(Y(\cdots Y(Y(a^{1},z_1)a^{2},z_2)\cdots,
a^{i-1}),z_{i-1})a^{i}$$ for some $a^{j}\in V^1$ and $k_j\in \Z$.
From  (\ref{eq-de-action}), (\ref{weak-comm1}),
(\ref{eq-adjoint}) and the skew symmetry we have the following
\begin{eqnarray*}
&&\langle q^\prime, Y_Q(v,z_0+z_i)Y_Q(u,z_i)q\rangle\\
&=&{\rm
Res}_{z_1}\cdots {\rm Res}_{z_{i-1}}z_1^{k_1}\cdots
z^{k_{i-1}}_{i-1}\\
&&\cdot\langle q^\prime, Y_Q(v,z_0+z_i)Y_Q(Y(\cdots
Y(Y(a^{1},z_1)a^{2},z_2)\cdots,
z_{i-1})a^{i},z_i)q\rangle\\
&=&{\rm Res}_{z_1}\cdots {\rm Res}_{z_{i-1}}z_1^{k_1}\cdots
z^{k_{i-1}}_{i-1}\langle
Y_Q^\prime(e^{(z_0+z_i)L(1)}\big(-(z_0+z_i)^{-2}\big)^{L(0)}v,(z_0+z_i)^{-1})q^\prime,\\
&& Y_Q(a^1,z_1+z_2+\cdots+z_i)\cdots Y_Q(a^{i-1},z_{i-1}+z_i)
Y_Q(a^i,z_i)q\rangle\\
&=&{\rm Res}_{z_1}\cdots {\rm Res}_{z_{i-1}}z_1^{k_1}\cdots
z^{k_{i-1}}_{i-1}\langle q^\prime,
Y_Q(a^1,z_1+z_2+\cdots+z_i)\cdots
 Y_Q(a^i,z_i) Y_Q(v,z_0+z_i)q\rangle\\
&=&{\rm Res}_{z_1}\cdots {\rm Res}_{z_{i-1}}z_1^{k_1}\cdots
z^{k_{i-1}}_{i-1}\\
&&\cdot\langle q^\prime, Y_Q(Y(Y(Y(\cdots
Y(Y(a^{1},z_1)a^{2},z_2)\cdots,
z_{i-1})a^{i},z_i-(z_0+z_i))v,z_0+z_i)q\rangle\\
&=&\langle q^\prime, Y_Q(Y(u,-z_0)v,z_0+z_i)q\rangle\\
&=&\langle q^\prime, Y_Q(e^{(-z_0)L(-1)}Y(v,z_0)u,z_0+z_i)q\rangle\\
&=&\langle q^\prime, Y_Q(Y(v,z_0)u,z_i)q\rangle.
\end{eqnarray*}
Thus we obtain the equation (\ref{eq-ass}).

Since $Q_i$ and  $Q_{i+1}$ are indecomposable $V^G$-modules, it
follows from   Corollary 3.10 of \cite{HLZ1} that for each $i$ there exists
$h_i\in\Q$ such that $\mathcal Y^i(u,z)Q_i\subseteq
z^{h_i}Q_{i+1}((z))$ for $u\in V^1$. Now the formula
(\ref{weak-comm1}) implies $h_i\equiv h_{i+1}\ ({\rm mod}\ \Z)$
and  the formula (\ref{commua-dia}) implies
$\sum_{i=0}^{p-1}h_i\equiv 0\ ({\rm mod}\ \Z)$. As a result, there
exists an integer $k$ $(0\leq k\leq p-1)$ such that $h_j\equiv
\frac{k}{p}\  ({\rm mod}\ \Z)$ for all $j$.  Thus $Q$ turns out to
be an admissible $g^l$-twisted $V$-module of finite length, where
$l$ is a positive integer such that $lk\equiv 1({\rm mod}\ \Z) $ if
$(p,k)=1$ and $l=(p,k)$ if $(p,k)\neq 1$. Clearly, $W$ is a $V^G$-homomorphic image of $Q$.

For the general case, since $V^G=(V^H)^{G/H}$ for any
subgroup $H$ of $G$, the proof can be carried out by induction on
$|G|$. So it suffices  to consider the case that $G$  is generated by two elements $g_1$
and $g_2$, say. Let $G_i$ be the subgroup of $G$ generated by $g_i$
of order $T_i$ for $i=1,2$.

  Let
$W$ be an indecomposable projective  $V^G$-module of finite
length. Then by the cyclic case we see that $W$ is  a
$V^G$-homomorphic image of an admissible  $g_2^i$-twisted
$V^{G_1}$-module of finite length for some $1\leq i\leq T_2$. So
by Lemma \ref{lemma-projective} we will be done if we can show
that each indecomposable projective  $g_2^i$-twisted
$V^{G_1}$-module of finite length is a homomorphic image of an
admissible $g_2^i$-twisted $V$-module. Here we only restrict ourself to
the case $i=1$, the other cases can be treated similarly.

Let $W$ be an indecomposable projective  $g_2$-twisted
$V^{G_1}$-module of finite length. Set $V^{(s,t)}=\{v\in
V|g_1v=e^{\frac{2s\pi i}{T_1}}v\ {\rm and}\ g_2v=e^{\frac{2t\pi
i}{T_2}}v\}$ for $1\leq s\leq T_1$, $1\leq t\leq T_2$ and
$V^{1}=\oplus_{j=1}^{T_2}V^{(1,j)}$. Then we have the following
$V^{G_1}$-epimorphism
$$\underbrace{V^1\boxtimes V^1\boxtimes\cdots \boxtimes V^1}_{T_1}\boxtimes W\rightarrow W,$$
from which the $V^{G_1}$-module $$Q=\bigoplus_{i=0}^{T_1-1}Q_i$$ can
be obtained similarly.  Let $\mathcal Y^i$ be the
canonical logarithmic intertwining operator of type
$\left(\begin{array}{cc}
    Q_{i+1} \\
  V^1\ Q_i \\
\end{array}\right)$ for $0\leq i\leq T_1-1$. By the similar argument as in the cyclic case, there
exist $h_{ij}\in\Q$ such that $\mathcal Y^i(u,z)Q_i\subseteq
z^{h_{ij}}Q_{i+1}((z))$ for  $u\in V^{(1,j)}$,
$0\leq i\leq T_1-1$  and $1\leq j\leq T_2$. Analogously, one can
show that $h_{iq}\equiv h_{jq}\ ({\rm mod}\ \Z)$ for
all $0\leq i,j\leq T_1-1$, $1\leq q\leq T_2$ and that
$\sum_{i=0}^{T_1-1}h_{il_i}\equiv \sum_{i=0}^{T_1-1}\frac{l_i}{T_2} \, ({\rm mod}\, \Z)$. It follows these two relations that  $\mathcal
Y^i(u^j,z)Q_i\subseteq z^{\frac{j}{T_2}}Q_{i+1}((z))$ for all
$u^j\in V^{(1,j)}$, $j\in\{1,2,\cdots, T_2\}$ and
$i\in\{0,1,\cdots,T_1-1\}$.  Since $V^1$ generates $V$, the action
of $V^{G_1}\oplus V^1$ on $Q$ can be extended to the whole $V$.
After equipping with such action, one can see that $Q$ carries the
structure of an admissible $g_2$-twisted $V$-module. This completes the proof.\QED
\bigskip

{\bf Proof of Theorem \ref{conjec}.}
 By  Lemma \ref{lemma-projective} and  Theorem
\ref{thm-submodule} we know that each irreducible $V^G$-module $L$
is a $V^G$-homomorphic image of  $\sum_{g\in G} M_g$, where $M_g$
is an admissible $g$-twisted $V$-module   for any $g\in G$. Since
$V$ is $g$-rational for any $g\in G$,  then it follows from
Theorem 2 of \cite{MT} that  $\sum_{g\in G} M_g$ is completely
reducible as $V^G$-module. So  $L$ is  a submodule of $\sum_{g\in
G} M_g$ and hence  of an irreducible $g$-twisted $V$-module for
some $g\in G$. We have proved Theorem \ref{conjec}. \QED

\bigskip
Let $G$ be a cyclic group of order $T$ generated by , say $g$, in
Theorem \ref{conjec}. Then   it is natural to ask that if it can
happen that there exists an irreducible $V^G$-module $L$ such that
 $L$ is not only contained in an irreducible $g^i$-twisted $V$-module
 but also in an irreducible $g^j$-twisted $V$-module for $0\leq i\neq j\leq T-1$. Set $V^1=\{v\in V\,|\, gv=e^{\frac{2\pi i}{T}}\}$.
 The following result tells us that this will not happen.

\begin{prop}
Let $V$ be a simple vertex operator superalgebra of CFT type such
that $V$ is $g$-rational for any $g\in G$ and  $V^G$ is
$C_2$-cofinite, where $G$ is a cyclic subgroup of $Aut(V)$ with
finite order $T$.  Let $L$ be an irreducible $V^G$-module.
Suppose that
$$\underbrace{V^1\boxtimes V^1\boxtimes\cdots \boxtimes V^1}_{T}=V^G.$$
Then there exists a  unique $i\ (0\leq i\leq T-1)$ such that  $L$
is contained in an irreducible $g^i$-twisted $V$-module $M^i$.
\end{prop}
\noindent{\em Proof.}\ \ By Theorem \ref{conjec} there exists $i\
(0\leq i\leq T-1)$ such that $L$ is a $V^G$-submodule of an
irreducible $g^i$-twisted $V$-module $M^i$.  Suppose that
$L\subset M^j$ for $0\leq i\neq j\leq T-1$. We first show that
$V^1\boxtimes L$ is an irreducible $V^G$-module. By the
associative property of tensor products and the assumption
$\underbrace{V^1\boxtimes V^1\boxtimes\cdots \boxtimes
V^1}_{T}=V^G$, we have $$\underbrace{V^1\boxtimes
(V^1\boxtimes\cdots \boxtimes(V^1}_{T}\boxtimes L)\cdots)\cong
V^G\boxtimes L\cong L.$$ This implies the  irreducibility of
$V^1\boxtimes L$, since any nonzero proper $V^G$-submodule $W$ of
$V^1\boxtimes L$ will give a nonzero proper submodule
$$\underbrace{V^1\boxtimes (V^1\boxtimes\cdots \boxtimes
(V^1}_{T-1}\boxtimes W)\cdots)$$ of $L$ (see Proposition
\ref{maintheorem-1} below).

We remain notation as in Lemma \ref{lemma-1}. Now it follows from
the universal property of the tensor product $V^1\boxtimes L$ and
the definition of $V^1\cdot L$ that there exists a
$V^G$-homomorphism $\phi_i$ from $V^1\boxtimes L$ onto $V^1\cdot
L\subset M^i$. In fact, $\phi_i$ is isomorphic, since
$V^1\boxtimes L$ is irreducible. Similarly, $V^1\boxtimes L\cong
V^1\cdot L\subset M^j$. Thus the $V^G$-submodule $\<u_nw\,|\,u\in
V^1, w\in L\ {\rm and}\ n\in\frac iT+\Z\>$ of $M^i$ and the
$V^G$-submodule $\<u_nw\,|\,u\in V^1, w\in L\ {\rm and}\
n\in\frac jT+\Z\>$ of $M^j$ are isomorphic, since both are
isomorphic to $V^1\boxtimes L$. On the other hand, notice that the
conformal weights of these two irreducible $V^G$-modules differ
$\frac{i-j}{T}$ from each other, so they can not be isomorphic as
$V^G$-modules, a contradiction. Thus $L$ is  contained only in
$g^i$-twisted $V$-modules. \QED
\bigskip

\section{Equivalence of rationalities}
\setcounter{section}{4}\setcounter{thm}{0}\setcounter{equation}{0}
\bigskip

This section is devoted to the proof of Theorem 1. We first have:

\begin{thm}\label{maintheorem-2}
Let $V$ be a simple vertex operator superalgebra and $G$ a finite
solvable subgroup of $Aut(V)$. Suppose that $V^G$ is rational.
Then $V$ is $g$-rational for any $g\in G$.
\end{thm}
\noindent{\em Proof.} Let $M$ be an admissible $g$-twisted
$V$-module. First,  consider the case that $G$ is abelian. Since
$V^G$ is rational, we can take an irreducible $V^G$-submodule $L$
of $M$. Consider the $V$-submodule $W=\sum_{\chi\in{\rm
Irr}G}V^{\chi}\cdot L$ of $M$. It is enough to show that $W$ is a
direct sum of irreducible $g$-twisted $V$-modules.

 Denote  $\tilde{L}$ the sum of all irreducible $V^G$-submodules
of $W$ isomorphic to $L$. Set $A=\{\chi\in{\rm Irr}G\,|\,L\cong
V^{\chi}\cdot L\}$. It is easy to see that $A$ is a subgroup of
Irr$G$ and $U=\sum_{\chi\in A}V^\chi$ is a vertex operator
subalgebra of $V$¡¡containing $V^G$. Clearly, the restriction of
$g$ to $U$ is an automorphism of $U$ and $\tilde L$ is a
$g$-twisted $U$-module. In fact, we are going to show that $\tilde
L$ is a completely reducible $g$-twisted $U$-module.

Note that $\tilde L$ can be  written   as ${\rm Hom}_{V^G}(L,
\tilde L)\otimes L$.¡¡ By the definition of $A$, $L\cong V^\chi
\cdot L$ for any $\chi\in A$.  Let $f_{\chi}$ be a
$V^G$-isomorphism from $L$ to $V^\chi \cdot L$ such that
$f_{1}={\rm Id}_L$. Then such $f_{\chi}$'s linearly span the
vector space ${\rm Hom}_{V^G}(L, \tilde L)$.

Next we shall extend each $f_\chi$ to a $U$-isomorphism $\hat
f_\chi$ of $\tilde L$. For this, define a liner map
$f_{\chi,\lambda}$ from $V^{\lambda}\cdot L$ to
$V^{\chi\lambda}\cdot L$ such that
$Y(a,z)f_\chi(w)=f_{\chi,\lambda}(Y(a,z)w)$ for $a\in V^\lambda$
and $w\in L$. Note that  for any $b\in V$ and $w\in L$ there
exists a nonnegative integer $l$ such that  the following weak
associativity
$$(z_0+z_2)^lY(Y(b,z_0)a,z_2)w=(z_0+z_2)^lY(b,z_0+z_2)Y(a,z_2)w$$ holds true for any $a\in V$.
By taking $b\in V^{\chi^{-1}}$ and $a\in V^\chi$ in the  formula
above one can verify  that the map $f_{\chi, \lambda}$ is
well-defined, that is, $a_n w=0$ implies $a_nf_{\chi}(w)=0\
(n\in\Q)$. As a result, $f_{\chi,\lambda}$ is not zero. Moreover,
using the weak associativity again one can show that
$f_{\chi,\lambda}$ is a $V^G$-homomorphism. Thus these two remarks
imply that $f_{\chi,\lambda}$ is a $V^G$-isomorphism, since by
Lemma \ref{lemma-1} both $V^\lambda\cdot L$ and
$V^{\chi\lambda}\cdot L$ are irreducible. For simplicity, we shall
take $f_{\chi,\lambda}=f_{\chi,\xi}$ whenever $V^\lambda\cdot
L=V^\xi\cdot L$. Define $\hat f_{\chi}: \tilde L\rightarrow \tilde
L$ by $\hat f_{\chi}(u)=f_{\chi,\lambda}(u)$ for $u\in
V^\lambda\cdot L$ and  $\lambda\in{\rm Irr}G$. It follows from a
similar argument as in the proof of Theorem 4.4 of \cite{DLM3}
 that $\hat f_\chi$ is a $U$-isomorphism  of $\tilde
L$.

Since both $f_{\lambda} f_{\chi}$ and $f_{\lambda\chi}$ are
$V^G$-isomorphisms from $L$ to $V^{\lambda\chi}$, there exists a
scalar $\beta(\lambda,\chi)\in \C^*$ such that
$f_{\lambda}f_{\chi}=\beta(\lambda, \chi)f_{\lambda\chi}$. One can
check that $\beta$ is a 2-cocycle. Hence we can define a
representation of $\C[A]_\beta$ on ${\rm Hom}_{V^G}(L, \tilde L)$
by sending $\chi $ to $f_\chi$, where $\C[A]_\beta$ is the twisted
group algebra associated with $\beta$. So $\tilde L={\rm
Hom}_{V^G}(L, \tilde L)\otimes L$ turns out to be a
$\C[A]_\beta\otimes V^G$-module. Then  from the argument  in the
proof of Theorem 3.10 of \cite{Lam}  we see that for any
$\C[A]_\beta$-submodule $Q$ of ${\rm Hom}_{V^G}(L, \tilde L)$,
$Q\otimes L$ is a $U$-module and that $Q\otimes L$ is
irreducible if and only if $Q$ is irreducible. It follows from the
complete reducibility of ${\rm Hom}_{V^G}(L, \tilde L)$ as
$\C[A]_\beta$-module, we can decompose $\tilde L$ into a direct
sum of irreducible $U$-modules, say $\tilde L=Q_1\otimes L\oplus
Q_2\otimes L\oplus\cdots\oplus Q_k\otimes L$.

Let $W_i=\sum_{\chi\in{\rm Irr}G}V^\chi\cdot (Q_i\otimes L)$ for $i=1,2,\cdots, k$. We
assert that each $W_i$ is an irreducible $g$-twisted $V$-module.
Let $X_i$ be a nonzero $V$-submodule of $W_i$, then $X_i\bigcap
V^\chi\cdot (Q_i\otimes L)\neq 0$ for some $\chi\in {\rm Irr}G$.
Applying the action of $V^{\chi^{-1}}$ on this sapce, we have
$$0\neq V^{\chi^{-1}}\cdot\big(X_i\bigcap V^\chi\cdot (Q_i\otimes L)\big)\subseteq X_i\bigcap Q_i\otimes L.$$
Now it follows from the irreducibility of $Q_i\otimes L$ that
$Q_i\otimes L\subseteq X_i$. As a result,  one has $W_i=X_i$, as $W_i$
is generated by $Q_i\otimes L$.  Whence $W=\sum_{i=1}^kW_i$ is a sum
of irreducible $g$-twisted $V$-modules. This completes the proof of
the theorem for the abelian case.

Since $G$ is solvable, there exists a normal subgroup $H$ of $G$ such that $G/H$ is abelian. Note that $V^G=(V^H)^{G/H}$. By the abelian case, $V^H$ is $g$-rational for any $g\in G$. Let $M$ be any admissible $g$-twisted $V$-module, which naturally becomes an
admissible $g$-twisted $V^H$-module. Then take an irreducible $V^H$-submodule $L$ of $M$ and  set $$W=\sum_{\chi\in\ {\rm Irr}(G/H)} (V^H)^\chi\cdot L.$$
 A similar
argument as in the abelian case shows that $W$ is completely reducible. Hence, $V$ is $g$-rational for any $g\in G$.  The proof is complete.\QED
\bigskip

\begin{prop}\label{maintheorem-1}
Let $V$ be a simple vertex operator superalgebra of CFT type and
$G$ a finite abelian subgroup of Aut$(V)$. Suppose that $V^G$ is
$C_2$-cofinite and $V$ is $g$-rational for any $g\in G$.  Then
$V^G$ is rational.
\end{prop}

\noindent{\em Proof.} Since $V^G$ is of CFT type,  the
$C_2$-cofiniteness of $V^G$ implies that any admissible
$V^G$-module $M=\oplus_{n\in\Z}M(n)$ generated by one element has
the following property:  each  subspace $M(n)$ is finite
dimensional (see Theorem 1 of \cite{B} or Lemma 2.4 of \cite{M1}).
Then by the argument in Lemma \ref{lemma-projective} we see that
$M\in \mathcal F_{V^G}$.

To prove the rationality of $V^G$, it is enough to show the
complete reducibility of $M$. But this follows immediately from
the proof of Theorem \ref{conjec}. \QED

\bigskip
 Combing Proposition \ref{maintheorem-1} and Theorem
\ref{maintheorem-2} gives Theorem \ref{t:I}.

\section{An interpretation of affine partition function of type $D_{2\rho+1}$}

\setcounter{section}{5}\setcounter{thm}{0}\setcounter{equation}{0}

 The classification of all possible affine
partition functions in terms of $A_1^{(1)}$ characters which are
$SL(2,\Z)$-invariant is given in Table 1 of \cite{CIZ}. They have
types
$$A_{n-1}, D_{2\rho+2}, D_{2\rho+1}, E_6,E_7,E_8.$$
All these partition functions can be realized as the partition
functions of the extensions of vertex operator algebras $L(k,0)$
for nonnegative integers $k$  or the corresponding conformal nets
\cite{KaL} except for types $D_{2\rho+1}$ and $E_7.$ In this
section we explain from the theory of vertex operator
superalgebras why the affine partition function
$$\sum\limits_{\lambda\, {\rm
odd}=1}^{4\rho-1}|\chi_{\lambda}(\tau)|^2+\sum\limits_{\lambda\, {\rm
even}=2}^{4\rho-2}\chi_{\lambda}(\tau)\chi_{4\rho-\lambda}(\tau)^*$$ of
type $D_{2\rho+1}$ is $SL(2,\Z)$-invariant.

Let $L(k,0)$ be the simple affine vertex operator algebra
associated with $A^{(1)}_1$ of level $k$. It is well known (cf.
\cite{FZ} and \cite{Li2}) that $L(k,0)$ is rational if and only if
$k\in\Z_{\geq0}$, and in this case $\{L(k, i)=\oplus_{n=0}^\infty
L(k,i)_{\lambda_i+n}|\ i=0,1,\cdots, k\}$ gives a complete list of
inequivalent irreducible $L(k,0)$-modules, where the conformal
weight $\lambda_i$ of $L(k,i)$ equals $\frac{i(i+2)}{4(k+2)}$.
Here  we need to point out  that  $L(k,k)$ is self-dual, i.e.,
$L(k,k)\cong L(k,k)^\prime$.

 From now on, we only consider the case $4|(k-2)$. Under such
circumstance, one will immediately see $L(k,k)=\oplus_{n=0}^\infty
L(k,k)_{\frac12+n}$. Consider the extension  $L(k,0)\oplus L(k,k)$
of $L(k,0)$  by a self-dual irreducible module $ L(k,k)$, which
carries the structure of  a vertex operator superalgebra
\cite{Li3}. By Theorem \ref{maintheorem-2}, $L(k,0)\oplus L(k,k)$
is rational and $\sigma$-rational, where $\sigma\in$
$Aut\big(L(k,0)\oplus L(k,k)\big)$ is defined by requiring
$\sigma|_{L(k,0)}={\rm Id}$ and $\sigma|_{L(k,k)}=-{\rm Id}$. And
 all irreducible $L(k,0)\oplus L(k,k)$-modules and
irreducible $\sigma$-twisted $L(k,0)\oplus L(k,k)$-modules can be
classified as follows (cf. \cite{Li3}):

\begin{prop}\label{irr-classification} Let $k$ be a positive integer such that $4|(k-2)$. Then

(1) $\{L(k,i)\oplus L(k,k-i)|0\leq i<\frac k2, i\ {\rm even}\}$ is a
complete list of irreducible $L(k,0)\oplus L(k,k)$-modules;

(2) $\{L(k,i)\oplus L(k,k-i)|1\leq i< \frac k2 , i\ {\rm odd}
\}\cup \{L(k,\frac k2), \sigma \circ L(k,\frac k2):=\big(L(k,\frac
k2), Y_{L(k,\frac k2)\big)}(\sigma\cdot,x)\}$ is a complete list
of irreducible $\sigma$-twisted $L(k,0)\oplus L(k,k)$-modules.
$($We refer the reader to \cite{DZ1} for the notation of
$\sigma\circ L(k,\frac k2)$$).$
\end{prop}

Denote $\chi_{i+1}(\tau)= Tr_{L(k,i)}q^{L(0)-\frac{c}{24}}$ the
$q$-character of $L(k,i)$ for $0\leq i\leq k$, where $q=e^{2\pi
i\tau}$. For convenience, denote $L(k,i)\oplus L(k,k-i)$ by $M_i$
for $0\leq i\,\, {\rm even}<\frac k2$ and $W_i$ for $0\leq i\,\,
{\rm odd}<\frac k2$, and $W_{\frac k2}=L(k,\frac k2)$. Define a
linear operator $\sigma$ on $M_i$ by decreeing
$$\sigma|_{L(k,i)}={\rm Id}_{L(k,i)}\ \ {\rm and}\ \ \sigma|_{L(k,k-i)}=-{\rm Id}_{L(k,k-i)}.$$

 Let us recall the modular transformation law  for
$\chi_i(\tau)$ from  \cite{CIZ}. Write  $k$ as $4\rho-2$, where
$\rho\geq 2$. Set $N=2(k+2)=8\rho$. For $-4\rho\leq i\leq -1$,
define
\begin{equation*}\label{def-modul}\chi_{\lambda}(\tau)=\chi_{\lambda+8\rho}(\tau)=-\chi_{-\lambda}(\tau).\end{equation*}  Then
\begin{eqnarray}\label{tran-T}
T(\chi_{\lambda}(\tau))=
\chi_{\lambda}(\tau+1)=e(\frac{\lambda^2}{16\rho}-\frac18)\chi_{\lambda}(\tau)
\end{eqnarray}
and \begin{eqnarray}\label{tran-S}
S(\chi_{\lambda}(\tau))&=&\chi_{\lambda}(-\frac1\tau)=\frac{-i}{\sqrt{8\rho}}\sum_{\lambda^\prime\in\Z/N\Z}e(\frac{\lambda\lambda^\prime}{8\rho})\chi_{\lambda^\prime}(\tau)\nonumber\\
&=&\frac{-i}{\sqrt{8\rho}}\sum_{\lambda^\prime=1}^{4\rho-1}
\big(e(\frac{\lambda\lambda^\prime}{8\rho})\chi_{\lambda^\prime}(\tau)+e(\frac{-\lambda\lambda^\prime}{8\rho})\chi_{-\lambda^\prime}(\tau)\big)\nonumber\\
&=&\frac{-i}{\sqrt{8\rho}}\sum_{\lambda^\prime=1}^{4\rho-1}
\big(e(\frac{\lambda\lambda^\prime}{8\rho})\chi_{\lambda^\prime}(\tau)-e(\frac{-\lambda\lambda^\prime}{8\rho})\chi_{\lambda^\prime}(\tau)\big)\\
&=&\frac{-i}{\sqrt{8\rho}}\sum_{\lambda^\prime=1}^{4\rho-1}
2i\sin\frac{2\pi\lambda\lambda^\prime}{8\rho}\chi_{\lambda^\prime}(\tau)\nonumber\\
&=&\frac{1}{\sqrt{2\rho}}\sum_{\lambda^\prime=1}^{4\rho-1}
\sin\frac{2\pi\lambda\lambda^\prime}{4\rho}\chi_{\lambda^\prime}(\tau)\nonumber,
\end{eqnarray}where $e(x)=e^{2\pi i x}$.

Computing the  trace and super trace functions associated with  irreducible
($\sigma$-twisted) modules listed in Proposition
\ref{irr-classification}, we have
\begin{eqnarray*}
&&\sum_{0\leq i\ {\rm even}<2\rho-1}(|Tr_{M_i}
q^{L(0)-\frac{c}{24}}|^2+|Tr_{M_i}\sigma
q^{L(0)-\frac{c}{24}}|^2)+\sum_{1\leq i\ {\rm odd}<
2\rho-1}|Tr_{W_i}q^{L(0)-\frac{c}{24}}|^2\nonumber\\
&&\, \, \,\,\, \, \,\, \,\, \, \, \, \,
\,\,\,\,\,\,\,\,\,\,\,\,\,\,\,
+|Tr_{W_{2\rho-1}}q^{L(0)-\frac{c}{24}}|^2+|Tr_{\sigma\circ
W_{2\rho-1}} q^{L(0)-\frac{c}{24}}|^2\\
&=&\sum_{\lambda=1}^{4\rho-1}|\chi_\lambda(\tau)|^2+\sum\limits_{\lambda\,
{\rm
odd}=1}^{4\rho-1}|\chi_{\lambda}(\tau)|^2+\sum\limits_{\lambda\,
{\rm
even}=2}^{4\rho-2}\chi_{\lambda}(\tau)\chi_{4\rho-\lambda}(\tau)^*.
 \end{eqnarray*}

\begin{prop}\label{prop-CCCCC}The function \begin{equation}\label{mod-i}\sum_{\lambda=1}^{4\rho-1}|\chi_i(\tau)|^2+ \sum\limits_{\lambda\, {\rm
odd}=1}^{4\rho-1}|\chi_{\lambda}(\tau)|^2+\sum\limits_{\lambda\,
{\rm
even}=2}^{4\rho-2}\chi_{\lambda}(\tau)\chi_{4\rho-\lambda}(\tau)^*\end{equation}
is $SL(2,\Z)$-invariant.
\end{prop}
\noindent{\em Proof.\ \ } By (\ref{tran-T}) one can easily verify
that (\ref{mod-i}) is  $T $-invariant. For any integer $\triangle$
satisfying $2-4\rho\leq \triangle\leq 8\rho-2$ we have
\begin{equation}\label{eq---A}
\sum_{\lambda=1}^{4\rho-1}\cos\frac{\pi\lambda\triangle}{4\rho}=\left\{\begin{array}{llll}
 -1\, \, \,\,\, &\triangle\ {\rm even}\neq 0, 4\rho
 \\[5pt]
  0 \, \, &\triangle\ {\rm odd}\\[5pt]
 4\rho-1 \, \, &\triangle =0\\[5pt]
  -1 \, \, &\triangle =4\rho,\\
\end{array}\right.
\end{equation}

\begin{equation}\label{eq---B} \sum_{\lambda\, {\rm
odd}=1}^{4\rho-1}\cos\frac{\pi\lambda\triangle}{4\rho}=\left\{\begin{array}{llll}
 0\, \, &\triangle\neq 0, 4\rho
 \\[4pt]
 2\rho \, \,  &\triangle=0\\[4pt]
 -2\rho \, \, &\triangle =4\rho
\end{array}\right.
\end{equation}
and \begin{equation}\label{eq---C} \sum_{\lambda\, {\rm
even}=2}^{4\rho-2}\cos\frac{\pi\lambda\triangle}{4\rho}=\left\{\begin{array}{llll}
 0 &&\triangle\ {\rm odd}
 \\
 -1   &&\triangle\ {\rm even}\neq 0, 4\rho\\
 2\rho-1  &&\triangle =0\,\, {\rm or}\,\, 4\rho.
\end{array}\right.
\end{equation}  Applying $S$ to
(\ref{mod-i}) and by using
(\ref{tran-S}), (\ref{eq---A})-(\ref{eq---C}) we obtain
\begin{eqnarray*}
&&S\big(\sum_{\lambda=1}^{4\rho-1}|\chi_{\lambda}|^2+\sum\limits_{\lambda\,
{\rm
odd}=1}^{4\rho-1}|\chi_{\lambda}(\tau)|^2+\sum\limits_{\lambda\,
{\rm
even}=2}^{4\rho-2}\chi_{\lambda}(\tau)\chi_{4\rho-\lambda}(\tau)^*\big)\\
&=&\sum_{\lambda=1}^{4\rho-1}\frac1{2\rho}\sum_{\mu,\gamma=1}^{4\rho-1}\sin{\frac{\pi\lambda\mu}
{4\rho}}\sin{\frac{\pi\lambda\gamma}{4\rho}}\chi_{\mu}(\tau)\chi_{\gamma}(\tau)^*+\sum_{\lambda\,
{\rm
odd}=1}^{4\rho-1}\frac{1}{2\rho}\sum_{\mu,\gamma=1}^{4\rho-1}\sin\frac{\pi\lambda\mu}{4\rho}\sin\frac{\pi\lambda\gamma}{4\rho}
\chi_{\mu}(\tau)\chi_{\gamma}(\tau)^*\\
&&+\sum_{\lambda\, {\rm
even}=2}^{4\rho-2}\frac{1}{2\rho}\sum_{\mu,\gamma=1}^{4\rho-1}\sin\frac{\pi\lambda\mu}{4\rho}\sin\frac{\pi(4\rho-\lambda)\gamma}{4\rho}
\chi_{\mu}(\tau)\chi_{\gamma}(\tau)^*\\
&=&\frac1{4\rho}\sum_{\mu,\gamma=1}^{4\rho-1}\left\{\sum_{\lambda=1}^{4\rho-1}
\big(\cos\frac{\pi\lambda(\mu-\gamma)}{4\rho}-\cos\frac{\pi\lambda(\mu+\gamma)}{4\rho}\big)\right.\\
&&+ \sum_{\lambda\,
{\rm odd}=1}^{4\rho-1}\big(\cos\frac{\pi\lambda(\mu-\gamma)}{4\rho}-\cos\frac{\pi\lambda(\mu+\gamma)}{4\rho}\big)\\
&&+\left.\sum_{\lambda\, {\rm
even}=2}^{4\rho-2}\big(\cos(\frac{\pi\lambda(\mu+\gamma)}{4\rho}-\pi\gamma)-\cos(\frac{\pi\lambda(\mu-\gamma)}
{4\rho}+\pi\gamma)\big)\right\}\chi_{\mu}(\tau)\chi_{\gamma}(\tau)^*\\
&=&\frac1{4\rho}\sum_{\mu,\gamma=1}^{4\rho-1}\left\{\sum_{\lambda=1}^{4\rho-1}
\big(\cos\frac{\pi\lambda(\mu-\gamma)}{4\rho}-\cos\frac{\pi\lambda(\mu+\gamma)}{4\rho}\big)\right.\\
&&+ \sum_{\lambda\,
{\rm odd}=1}^{4\rho-1}\big(\cos\frac{\pi\lambda(\mu-\gamma)}{4\rho}-\cos\frac{\pi\lambda(\mu+\gamma)}{4\rho}\big)\\
&&+\left.(-1)^\gamma\sum_{\lambda\, {\rm
even}=2}^{4\rho-2}\big(\cos(\frac{\pi\lambda(\mu+\gamma)}{4\rho})-\cos(\frac{\pi\lambda(\mu-\gamma)}
{4\rho})\big)\right\}\chi_{\mu}(\tau)\chi_{\gamma}(\tau)^*\\
&=&\sum_{\lambda=1}^{4\rho-1}|\chi_{\lambda}|^2+\sum\limits_{\lambda\,
{\rm
odd}=1}^{4\rho-1}|\chi_{\lambda}(\tau)|^2+\sum\limits_{\lambda\,
{\rm
even}=2}^{4\rho-2}\chi_{\lambda}(\tau)\chi_{4\rho-\lambda}(\tau)^*.
\end{eqnarray*}
That is,
$$\sum_{\lambda=1}^{4\rho-1}|\chi_{\lambda}|^2+\sum\limits_{\lambda\,
{\rm
odd}=1}^{4\rho-1}|\chi_{\lambda}(\tau)|^2+\sum\limits_{\lambda\,
{\rm
even}=2}^{4\rho-2}\chi_{\lambda}(\tau)\chi_{4\rho-\lambda}(\tau)^*$$
is $S$-invariant. Thus, it is $SL(2,\Z)$-invariant. \QED

\begin{rem} Proposition \ref{prop-CCCCC} is the main motivation for  Conjecture \ref{conje-AAA}. Although the modular invariance of trace functions for vertex operator superalgebras are known \cite{DZ2}, but there is still a distance to establish  the conjecture. As in \cite{DLN} one needs to prove that the $S$-matrix on the space spanned by
$$Tr_{M_{i}}q^{L(0)-\frac{c}{24}},
Tr_{M_{i}}\sigma q^{L(0)-\frac{c}{24}}, Tr_{W_{j}}q^{L(0)-\frac{c}{24}}$$
for all $i,j$ (see Introduction for notation $M_i,W_j$) is unitary.
\end{rem}

\begin{coro} The function $$\sum\limits_{\lambda\, {\rm
odd}=1}^{4\rho-1}|\chi_{\lambda}(\tau)|^2+\sum\limits_{\lambda\,
{\rm
even}=2}^{4\rho-2}\chi_{\lambda}(\tau)\chi_{4\rho-\lambda}(\tau)^*$$
is $SL(2,\Z)$-invariant.
\end{coro}

{\bf Proof.} In the  calculations above, we have already shown that \begin{equation}\label{eq--Las}\sum_{i=1}^{4\rho-1}|\chi_i(\tau)|^2\end{equation} is
 $SL(2,\Z)$-invariant (see also Corollary 3.8 of \cite{DLN}). Note that $$\sum\limits_{\lambda\, {\rm
odd}=1}^{4\rho-1}|\chi_{\lambda}(\tau)|^2+\sum\limits_{\lambda\,
{\rm
even}=2}^{4\rho-2}\chi_{\lambda}(\tau)\chi_{4\rho-\lambda}(\tau)^*$$
is the difference of the functions (\ref{mod-i}) and (\ref{eq--Las}).
The result follows from Proposition \ref{prop-CCCCC} immediately. \QED

\end{document}